\documentclass[10pt]{amsart}
\usepackage{amsmath}

\newtheorem{thm}{Theorem}
\newtheorem{lem}[thm]{Lemma}

\theoremstyle{definition}
\newtheorem*{remark}{Remark}

\newcommand{\Z}{{\mathbb Z}}

\newcommand{\N}{{\mathbb N}}
\newcommand{\T}{{\mathbb T}}

\newcommand{\supp}{\operatorname{supp}}
\renewcommand{\mod}{\;\operatorname{mod}}
\newcommand{\Tr}{\operatorname{Tr}}
\newcommand{\Det}{\operatorname{Det}}

\newcommand{\intinf}{\int_{-\infty}^\infty}

\newcommand{\E}{{\mathbb E}} 

\newcommand{\I}{1\!\!1} 

\renewcommand{\i}{{\mathrm{i}}} 
\renewcommand{\d}{{\mathrm{d}}} 
\renewcommand{\^}{\widehat}

\newcommand{\Zf}{Z_f} 

\newcommand{\U}{{\mathrm{U}}} 
\newcommand{\SO}{{\mathrm{SO}}} 
\newcommand{\SP}{\mathrm{Sp}} 

\renewcommand{\t}{\theta}

\begin{document}
\title[Mock--Gaussian Behaviour in Classical Compact Groups]{Mock--Gaussian Behaviour for Linear Statistics of Classical Compact Groups}
\author{C.P. Hughes and Z. Rudnick}
\address{Raymond and Beverly Sackler School of Mathematical Sciences,
Tel Aviv University, Tel Aviv 69978, Israel
({\tt hughes@post.tau.ac.il})}
\address{Raymond and Beverly Sackler School of Mathematical Sciences,
Tel Aviv University, Tel Aviv 69978, Israel
({\tt rudnick@post.tau.ac.il})}

\date{27 June 2002}

\thanks{Supported in part by  the EC TMR network 
"Mathematical aspects of Quantum Chaos", EC-contract no
HPRN-CT-2000-00103}

\begin{abstract}
We consider  the scaling limit of 
linear statistics for eigenphases of a matrix taken from one of the
classical compact groups.   We compute their moments 
and find that the first few moments are Gaussian, whereas the
limiting distribution is not. 
The precise number of Gaussian moments depends upon the particular
statistic considered. 
\end{abstract}

\maketitle
\section{Introduction}

In this paper we investigate the scaling limit of linear statistics
for eigenphases of matrices in the classical groups.  Given a unitary
$N\times N$ matrix $U$ with eigenvalues $e^{\i\theta_n}$, $1\leq n\leq
N$, and a test function $g$ which we assume is $2\pi$--periodic, 
consider the linear statistic
$$
\Tr g(U) := \sum_{n=1}^N g(\theta_n)
$$
A number of authors have studied the limiting distribution as $N\to
\infty$ of $\Tr g(U)$ as $U$ 
varies over a family $G(N)$ of classical groups and have concluded
that the distribution is Gaussian, see \cite{DS, DE, Johan}. 

Soshnikov \cite{Sosh} showed that this result remains valid in the
``mesoscopic'' regime, that is if one 
considers eigenphases $\theta_n$ in an interval of length about
$1/L$ where $L=L_N\to \infty$ but $L/N\to 0$:  
For a Schwartz function $f$ on the real line, define 
$$
F_L(\theta):=\sum_{j=-\infty}^\infty f(\frac L{2\pi}(\theta+2\pi j))
$$
which is $2\pi$-periodic and localised on a scale of $1/L$. 
Soshnikov \cite{Sosh}
showed that as long as $L/N\to 0$, then the limiting distribution of
$\Tr F_L(U)$ as $U$ ranges over all unitary matrices in $\U(N)$, 
$N\to\infty$  
is a Gaussian with mean 
$$
\frac{N}{L} \int_{-\infty}^{\infty}f(x)\;\d x
$$
and variance  
$$
\intinf \^ f(t)^2 |t|\d t
$$
where the Fourier transform is defined as
$$
\^ f(t) := \intinf f(x) e^{-2\pi\i xt}\;\d x
$$

There are similar formulae for the other classical groups. 

Our goal is to investigate these linear statistics in the {\em scaling
limit}, that is to take $L=N$. Thus we set 
$$
\Zf(U):=\Tr F_N(U) = \sum_{n=1}^N F_N(\theta_n)
$$
In \cite{hr1} we proved
\begin{thm}\label{thm:mock_Gauss_unitary}
If  $\supp \^f \subseteq [-2/m,2/m]$   then the first
$m$ moments of $\Zf(U)$ over the unitary group $\U(N)$ converge as
$N\to \infty$ 
to the Gaussian moments with mean
$\intinf f(x)\d x$ and variance 
$$ \int_{-\infty}^\infty \min(|u|,1) |\^f(u)|^2 \;\d u$$  
\end{thm}
We called this a ``mock-Gaussian'' behaviour. It is worth remarking
that in \cite{hr1} we find the full distribution of $\Zf$, and it is
not Gaussian. Only the first few moments are.

The purpose of this paper is to demonstrate mock-Gaussian behaviour for
linear statistics in other classical compact groups, the special
orthogonal group $\SO(N)$ and
the symplectic group $\SP(N)$ ($N$ must be even in the symplectic
group). If $e^{\i\t}$ is an eigenvalue of a matrix $U$ taken from   
one of these groups then $e^{-\i\t}$ is an eigenvalue too. This means $1$ is
always an eigenvalue of $U\in\SO(N)$ if $N$ is odd.

Due to the pairing of eigenvalues, the function $f$ must be even. 
Our results are
\begin{thm}\label{thm:mock_Gaussian}

i) If $\supp \^f \subseteq [-1/m,1/m]$ then the first
$m$ moments of $\Zf(U)$ over the symplectic group $\SP(N)$ 
converge to the Gaussian moments with mean  
$$\^f(0)-\int_{0}^1 \^f(u)\;\d u$$
and variance
$$2\int_{-1/2}^{1/2} |u| \^f(u)^2\;\d u$$

ii) If $\supp \^f \subseteq (-1/m,1/m)$ then the first $m$ moments of
$\Zf(U)$ over the special orthogonal group  $U\in\SO(N)$ 
converge to the Gaussian moments with mean  
$$\^f(0)+\int_{0}^1 \^f(u)\;\d u$$
and variance
$$2\int_{-1/2}^{1/2} |u| \^f(u)^2\;\d u$$ 

\end{thm}

\begin{remark}
There exists $f$ such that $\supp \^f \subseteq [-1/m,1/m]$ and whose
$m+1$-st moment is not Gaussian.
\end{remark}

\subsection{Moments and cumulants}
One approach to proving such results is to use the Fourier expansion
$g(\theta) = \sum_n g_n e^{\i n \theta}$ and expand $\Tr g(U)$ as a sum
$$
\Tr g(U) = \sum_n g_n \Tr(U^n)
$$
Computing moments of $\Tr g(U)$ then boils down to being able to compute
integrals of products of $\Tr(U^n)$ over the classical group. 
Theorem~\ref{thm:mock_Gauss_unitary} for the unitary group 
was proven in \cite{hr1} using this approach by employing a result of 
Diaconis and Shahshahani,
\cite{DS,DE}, concerning moments of traces of random unitary
matrices. Their result is a consequence of Schur duality for
representations of the unitary group and the symmetric group, and the
second orthogonality relation for characters of the symmetric group. 

The paper by Diaconis and Evans \cite{DE} (see also \cite{DS}) contains a corresponding result for moments of
traces of random symplectic and orthogonal\footnote{Note that Diaconis and Evans consider orthogonal matrices, whereas we are interested in the special orthogonal group} matrices (which they deduce using the work of Ram \cite{Ram} on Brauer algebras), which can be used
to prove our theorems in {\em half} the range, that is the $m$-th
moment of $\Zf$ is Gaussian if  $\supp \^f$ lies in the interval
$(-1/2m,1/2m)$.  
We wish to have the {\em full} range so as to compare with zeros of quadratic
$L$--functions, where linear statistics  show 
mock-Gaussian behaviour in the same  full range (this can be deduced from
the work of Rubinstein, \cite{Rubin}). 
The case of Dirichlet $L$--functions, which correspond to the unitary group, was considered in \cite{hr1}.

To obtain the results we desire, we abandon moments and instead use 
the {\em cumulants} $C_{\ell}^{G(N)}(g)$ of $\Tr g(U)$. 
These are defined via the expansion 
$$
\log \E_{G(N)}(e^{t\Tr g(U)}) = \sum_{\ell=1}^\infty C_{\ell}^{G(N)}(g) \frac{t^\ell}{\ell!}
$$ 
where $\E_{G(N)}$ denotes the expectation with respect to Haar measure
over the group $G(N)$. The cumulants have previously been considered in this context 
by Soshnikov \cite{Sosh} (interestingly, his results again only give 
half the required range), and it is his combinatorial approach that we adopt.

There is a
natural decomposition for the cumulants on the symplectic and special
orthogonal groups. For brevity we will describe the situation for the
symplectic group (so $N$, the matrix size, is assumed to be even). The cumulants can be written as
$$
C_\ell^{\SP(N)}(g) = 2^\ell C_{\ell,N+1}^{even}(g) - 2^\ell C_{\ell,N+1}^{odd}(g) 
$$

We show that the odd parts $C_{\ell,N+1}^{odd}(g)$    
of the cumulants vanish in a certain region, and in fact if $g_k=0$ for $|k|>(N+1)/\ell$ then the $\ell$-th cumulant vanishes.

For all $g$, 
the even summand equals a unitary cumulant: 
$$ 
C_{\ell,N+1}^{even}(g) = \frac 12 C_{\ell}^{\U(N+1)}(g)
$$
We may now employ  the available results about the unitary group to
deduce that $C_{\ell,N+1}^{even}(g)$ also vanishes in a larger region.
Setting $g=F_N$ we obtain Theorem~\ref{thm:mock_Gaussian}.

Since  moments and cumulants give essentially equivalent information,
we can now go back to computing averages of the product of traces on
classical groups and resolve a problem raised in \cite[Remark 8.2]{DE},
to show 

\begin{thm}\label{thm:mmts traces}
Let $Z_j$ be independent standard normal random variables, and let
\begin{equation*}
\eta_j = 
\begin{cases}
1 & \text{ if $j$ is even}\\
0 & \text{ if $j$ is odd}
\end{cases}
\end{equation*}

i) If $a_j\in\{0,1,2,\dots\}$ for $j=1,2,\dots$ are such that $\sum j a_j \leq N+1$, where $N$ is even, then
\begin{equation*}
\E_{\SP(N)}\left\{ \prod (\Tr U^j)^{a_j}\right\} = E\left\{\prod (\sqrt{j} Z_j-\eta_j)^{a_j}\right\}
\end{equation*}

ii) If $a_j\in\{0,1,2,\dots\}$ for $j=1,2,\dots$ are such that $\sum j a_j \leq N-1$ then 
\begin{equation*}
\E_{\SO(N)}\left\{ \prod (\Tr U^j)^{a_j}\right\} = E\left\{\prod (\sqrt{j} Z_j+\eta_j)^{a_j}\right\}
\end{equation*}
\end{thm}
Similar Theorems have been proven by Diaconis and Evans \cite{DE}, though only for half the
range (that is, they require $\sum ja_j \leq N/2$).

\section{Cumulants of linear statistics}

In order to calculate $C_\ell^{\SP(N)}(g)$ we need to know the moment
generating function. Weyl \cite{Weyl} showed that $\E_{\SP(N)} \{ e^{t \Tr g(U)} \}$ could be written as an integral over the $N/2$ independent eigenphases (recall that $N$ must be even for a symplectic matrix to exist). He showed that, writing $N=2M$,
\begin{align*}
\E_{\SP(N)} \{ e^{t \Tr g(U)} \} &= \E_{\SP(N)} \left\{ \exp\left(2t \sum_{n=1}^M g(\t_n) \right) \right\}\\
&= \int_{[0,\pi]^M} \Det\{ Q^{\SP(2M)}(\t_i,\t_j)\}_{1\leq i,j\leq M}\prod_{n=1}^M e^{2t g(\t_n)} \;\d\t_n
\end{align*}
where the kernel is $Q^{\SP(N)}(x,y) := S_{N+1}(x-y)-S_{N+1}(x+y)$ with
\begin{equation}\label{eq:defn S_N}
S_{N}(z) := \frac{1}{2\pi}\frac{\sin(Nz/2)}{\sin(z/2)}
\end{equation}

Now, it is a general fact that if $\t_n\in \T$, where $\T$ is some real interval, are such that
$$\E\left\{ \exp\left(\sum_{n=1}^M t g(\t_n)\right) \right\} = \int_{\T^M} \Det\{ Q_M(\t_i,\t_j)\}_{1\leq i,j\leq M}\prod_{n=1}^M e^{t g(\t_n)} \;\d\t_n$$
then if the $\ell$th cumulant of $\sum g(\t_n)$, $C_\ell$, is defined by the expansion
$$\log \E \left\{ \exp\left( t \sum_{n=1}^M g(\t_n)\right)\right\} = \sum_{\ell=1}^\infty \frac{t^\ell}{\ell!} C_{\ell}$$
then \cite{Sosh,Stojanovic}
\begin{equation}\label{eq:general cumulants formula}
C_{\ell} = \sum_{m=1}^\ell  \sum_{\sigma\in P(\ell,m)} (-1)^{m+1} (m-1)! \int_{\T^m} \prod_{j=1}^m g^{\lambda_j}(x_j) Q_M(x_j,x_{j+1}) \;\d x_j
\end{equation}
where we identify $x_{m+1}$ with $x_1$. Here $P(\ell,m)$ is the set of
all partitions of $\ell$ objects into $m$ non-empty blocks, where the $j$th block has
$\lambda_j=\lambda_j(\sigma)$ elements (that is $\lambda_j := \#\{i \ : \ 1\leq i\leq \ell \ , \ \sigma(i)=j\}$).

Thus,
\begin{multline*}
C_{\ell}^{\SP(N)}(g) =\\
2^\ell\sum_{m=1}^\ell  \sum_{\sigma\in P(\ell,m)} (-1)^{m+1} (m-1)! \int_{[0,\pi]^m} \prod_{j=1}^m g^{\lambda_j}(x_j) Q^{\SP(N)}(x_j,x_{j+1}) \;\d x_j
\end{multline*}
Since $Q^{\SP(N)}(x,y)$ is odd in both variables, $\prod_{j=1}^m Q^{\SP(N)}(x_j,x_{j+1})$ is even in all variables, and so, since $g$ is an even function, we may extend the integral to be over $[-\pi,\pi]$ and thus
\begin{multline*}
C_{\ell}^{\SP(N)}(g) = 2^\ell \sum_{m=1}^\ell \sum_{\sigma\in P(\ell,m)} (-1)^{m+1} (m-1)! \times\\
\times \frac{1}{2^m} \int_{[-\pi,\pi]^m} \prod_{j=1}^m g^{\lambda_j}(x_j) \left(S_{N+1}(x_j-x_{j+1})-S_{N+1}(x_j+x_{j+1})\right) \; \d x_j
\end{multline*}
and on expanding out the middle product on the bottom line,
\begin{align*}
C_{\ell}^{\SP(N)}(g) &= 2^\ell\sum_{m=1}^\ell \sum_{\sigma\in P(\ell,m)} (-1)^{m+1} (m-1)! \frac{1}{2^m} \sum_{\epsilon_1=\pm 1,\dots,\epsilon_m=\pm 1} \times\\
&\qquad\qquad\times \int_{[-\pi,\pi]^m} \prod_{j=1}^m g^{\lambda_j}(x_j)  \epsilon_j S_{N+1}(x_j-\epsilon_j x_{j+1}) \;\d x_j\\
&= 2^\ell C_{\ell,N+1}^{\mathrm{even}}(g) - 2^\ell C_{\ell,N+1}^{\mathrm{odd}}(g)
\end{align*}
where $C_{\ell,N+1}^{\mathrm{even}}(g)$ contains those terms with $\prod_{j=1}^m \epsilon_j = +1$ and $C_{\ell,N+1}^{\mathrm{odd}}(g)$ contains those terms with $\prod_{j=1}^m \epsilon_j = -1$.

Similarly one can calculate the other groups, using Weyl's calculation
of Haar measure, which is summarised in table \ref{table:kernels}.

\begin{table}[hbt]
\centering
\begin{tabular}{||c|c|c|c||}
\hline\hline
{}&{}&{}&{}\\
Group & $\Tr g(U)$ & Kernel $Q_M(x,y)$ & Range $\T$ \\
{}&{}&{}&{}\\
\hline\hline
{}&{}&{}&{}\\
$\U(N)$ & $\sum_{n=1}^N g(\t_n)$ & $S_N(x,y)$ & $(-\pi,\pi]$\\
{}&{}&{}&{}\\
{}&{}&{}&{}\\
$\SP(N)$ & $2\sum_{n=1}^M g(\t_n)$ & $S_{N+1}(x-y)-S_{N+1}(x+y)$ & $[0,\pi]$\\
$N=2M$ &{}&{}&{}\\
{}&{}&{}&{}\\
$\SO(N)$ & $2\sum_{n=1}^M g(\t_n)$ & $S_{N-1}(x-y)+S_{N-1}(x+y)$ & $[0,\pi]$\\
$N=2M$ &{}&{}&{}\\
{}&{}&{}&{}\\
$\SO(N)$ & $g(0)+2\sum_{n=1}^M g(\t_n)$ & $S_{N-1}(x-y)-S_{N-1}(x+y)$ & $[0,\pi]$\\
$N=2M+1$ &{}&{}&{}\\
{}&{}&{}&{}\\
\hline\hline
\end{tabular}
\caption{Kernels for Haar measure over the classical compact groups} \label{table:kernels}
\end{table}

\subsection{Summary}

Put
\begin{multline}\label{eq:C^even}
C_{\ell,M}^{\mathrm{even}}(g) = \sum_{m=1}^\ell \sum_{\sigma\in P(\ell,m)} (-1)^{m+1} (m-1)! \frac{1}{2^m} \sum_{\substack{\epsilon_1=\pm 1,\dots,\epsilon_m=\pm 1\\\prod\epsilon_j=+1}} \times\\
\int_{[-\pi,\pi]^m} \prod_{j=1}^m g^{\lambda_j}(x_j)  S_{M}(x_j-\epsilon_j x_{j+1}) \;\d x_j
\end{multline}
and
\begin{multline}\label{eq:C^odd}
C_{\ell,M}^{\mathrm{odd}}(g) = \sum_{m=1}^\ell \sum_{\sigma\in P(\ell,m)} (-1)^{m+1} (m-1)! \frac{1}{2^m} \sum_{\substack{\epsilon_1=\pm 1,\dots,\epsilon_m=\pm 1\\\prod\epsilon_j=-1}} \times\\
\int_{[-\pi,\pi]^m} \prod_{j=1}^m g^{\lambda_j}(x_j)  S_{M}(x_j-\epsilon_j x_{j+1}) \;\d x_j
\end{multline}
with $S_M$ defined in \eqref{eq:defn S_N}.

\begin{itemize}
\item For all $\ell$,
\begin{equation*}
C_{\ell}^{\SP(2M)}(g) = 2^{\ell} C_{\ell,2M+1}^{\mathrm{even}}(g) - 2^\ell C_{\ell,2M+1}^{\mathrm{odd}}(g)
\end{equation*}

\item For all $\ell$,
\begin{equation*}
C_{\ell}^{\SO(2M)}(g) = 2^{\ell} C_{\ell,2M-1}^{\mathrm{even}}(g) + 2^\ell C_{\ell,2M-1}^{\mathrm{odd}}(g)
\end{equation*}

\item For $\ell=1$,
\begin{equation*}
C_{1}^{\SO(2M+1)}(g) = 2 C_{1,2M}^{\mathrm{even}}(g) - 2 C_{1,2M}^{\mathrm{odd}}(g) + \sum_{k=-\infty}^\infty g_k 
\end{equation*}
and for all $\ell\geq 2$,
\begin{equation*}
C_{\ell}^{\SO(2M+1)}(g) = 2^{\ell} C_{\ell,2M}^{\mathrm{even}}(g) - 2^\ell C_{\ell,2M}^{\mathrm{odd}}(g)
\end{equation*}
\end{itemize}

In the next section, we will show that $C_{\ell,M}^{\mathrm{even}}(g)
= \frac{1}{2} C_{\ell}^{\U(M)}(g)$, and then we will calculate
$C_{\ell,M}^{\mathrm{odd}}(g)$, first in the case when $M$ is odd, and then
in the case when $M$ is even.

The results will show that
\begin{equation}\label{eq:defn mu}
C_{\ell}^{G(N)}(g) = \sum_{k\in \Z^\ell} \mu_{\ell}^{G(N)}(k_1,\dots,k_\ell) \prod_{j=1}^\ell g_{k_j}
\end{equation}
where $\mu_{\ell}^{G(N)}(k_1,\dots,k_\ell)$ is invariant under permutations of its arguments.

Combining the results from the next section proves the following theorems:

\begin{thm}\label{thm:cumulants Sp(2M)}
$$C_1^{\SP(2M)}(g) = 2M g_0 - 2\sum_{n=1}^{M} g_{2n}$$
$$C_2^{\SP(2M)}(g) = 4\sum_{n=1}^\infty \min(n,2M+1) g_n^2 - 4\sum_{k=M+1}^\infty g_k^2 - 8\sum_{l=1}^{M}\sum_{k=M+1}^\infty g_{k+l} g_{k-l}$$
and for $\ell\geq 3$, $\mu_{\ell}^{\SP(N)}(k_1,\dots,k_\ell) = 0$ if $\sum_{j=1}^\ell |k_j| \leq N+1$.
\end{thm}

\begin{thm}\label{thm:cumulants SO(N)}
When averaged over the special orthogonal group, the mean of $\Tr g(U)$ is
\begin{gather*}
C_1^{\SO(2M)}(g) = 2M g_0 + 2\sum_{n=1}^{M-1} g_{2n}\\
C_1^{\SO(2M+1)}(g) = (2M+1)g_0+2\sum_{n=1}^M g_{2n} + 2\sum_{n=2M+1}^\infty g_n 
\end{gather*}
and the variance is
\begin{gather*}
C_2^{\SO(2M)}(g) = 4\sum_{n=1}^\infty \min(n,2M-1) g_n^2 + 4\sum_{k=M}^\infty g_k^2 + 8\sum_{l=1}^{M-1}\sum_{k=M}^\infty g_{k+l} g_{k-l}\\
C_2^{\SO(2M+1)}(g) = 4\sum_{n=1}^\infty \min(n,2M) g_n^2 - 8\sum_{\substack{n=1\\n \text{ odd}}}^{2M-1} \sum_{\substack{m=2M+1\\m \text{ odd}}}^\infty g_{(m+n)/2} g_{(m-n)/2}
\end{gather*}
For $\ell\geq 3$, $\mu_{\ell}^{\SO(N)}(k_1,\dots,k_\ell) = 0$ if $\sum_{j=1}^\ell |k_j| \leq N-1$.
\end{thm}

\section{The combinatorial calculations}

\subsection{The calculation of $C_{\ell,M}^{\mathrm{even}}(g)$}

The following lemma was stated by Soshnikov in \cite{Sosh}:
\begin{lem}\label{lem:SO(2N) even}
For all $\ell$,
\begin{equation*}
C_{\ell,M}^{\mathrm{even}}(g) = \frac{1}{2} C_{\ell}^{U(M)}(g)
\end{equation*}
\end{lem}

\textbf{Proof.}
Symbolically denote
\begin{equation}\label{eq:integral for SO(2N)_even}
 \int_{[-\pi,\pi]^m} \prod_{j=1}^m g^{\lambda_j}(x_j) S_{M}(x_j-\epsilon_j x_{j+1}) \;\d x_j
\end{equation}
by $(\epsilon_1,\epsilon_2,\dots,\epsilon_m)$. If $\epsilon_1=1$ do nothing, but if $\epsilon_1=-1$ then
change variables to $x_2\mapsto -x_2$, and note that since $g$ and $S_{M}$ are even functions, and the integral over $x_2$ is over $[-\pi,\pi]$, then \eqref{eq:integral for SO(2N)_even} becomes $(+1,-\epsilon_2,\epsilon_3,\dots,\epsilon_m)$.

Observe that this achieves the following: If the initial situation was
$(-1,-1,\dots)$ then it becomes $(+1,+1,\dots)$ while if it was
$(-1,+1,\dots)$ it becomes $(+1,-1,\dots)$. Therefore there is either
the same number of $-1$'s in the set of $\epsilon$ or there are two
less $-1$'s. 

Now repeat for the new $\epsilon_2$, changing variables only if it is $-1$, and so on all the way up to $\epsilon_m$. Each time the action either leaves the number of $-1$'s unchanged or reduces it by 2. Since we started with an even number of $-1$'s in the set of $\epsilon$ this algorithm will terminate with \eqref{eq:integral for SO(2N)_even} equaling $(+1,+1,\dots,+1)$, which is independent of $\epsilon$. There are $2^{m-1}$ possible $\epsilon$ with an even number of $-1$'s, and so
\begin{multline*}
C_{\ell,M}^{\mathrm{even}}(g) = \sum_{m=1}^\ell \sum_{\sigma\in P(\ell,m)} (-1)^{m+1} (m-1)! \frac{1}{2^m} \times\\
\times  2^{m-1} \int_{[-\pi,\pi]^m} \prod_{j=1}^m g^{\lambda_j}(x_j) S_{M}(x_j- x_{j+1}) \;\d x_j
\end{multline*}
which we recognise as $\frac{1}{2} C_{\ell}^{U(M)}(g)$.
\qed

The cumulants of a random unitary matrix have previously been
calculated, essentially by Soshnikov \cite{Sosh}, but they can also be
deduced from the work of Diaconis and Shahshahani \cite{DS} 
and of Diaconis and Evans \cite{DE}.
\begin{thm}{\bf (Soshnikov).}\label{thm:Sosh_unitary}
Let $C_\ell^{\U(N)}$ be the $\ell$th cumulant of $\Tr g(U)$, averaged over all $N\times N$ unitary matrices with Haar measure. Then
$$C_1^{\U(N)}=N g_0$$
$$C_2^{\U(N)} = \sum_{\substack{n=-\infty\\n\neq 0}}^\infty \min(|n|,N) g_n g_{-n}$$
and for $\ell\geq 3$,
\begin{equation*}
\left|C_\ell^{\U(N)}(g)\right| \leq \text{const}_\ell \sum_{\substack{k_1+\dots+k_\ell=0\\|k_1|+\dots+|k_\ell|>2N}} |k_1| |g_{k_1}| \dots |g_{k_\ell}|
\end{equation*}
\end{thm}

\begin{remark}
The heart of the proof of this theorem is a deep combinatorial fact called the Hunt-Dyson formula.
\end{remark}

\begin{remark}
Actually, the error term in \cite{Sosh} has the sum running over all $k_1+\dots+k_\ell=0$ such that $|k_1|+\dots+|k_\ell|>N$.
But it is clear from equation 2.9 of \cite{Sosh} that there is no contribution to $C_\ell^{\U(N)}$ for $\ell\geq 3$ if $\sum k_i \I_{\{k_i>0\}} \leq N$ and if  $\sum -k_i \I_{\{k_i<0\}} \leq N$. Since the $k_i$ sum to zero, it must be that the sum over positive terms equals the sum over negative terms, and so this is the same as the condition that $\sum |k_i| \leq 2N$, as we have it in the theorem.
\end{remark}

\subsection{The calculation of $C_{\ell,2M+1}^{\mathrm{odd}}(g)$}

Observe from \eqref{eq:defn S_N} that
\begin{equation}\label{eq:fourier expansion of S_(2M+1) kernel}
S_{2M+1}(z)=\frac{1}{2\pi} \sum_{n=-M}^{M} e^{-\i nz} 
\end{equation}

\begin{lem}
One can calculate $C_{1,2M+1}^{\mathrm{odd}}(g)$ and $C_{2,2M+1}^{\mathrm{odd}}(g)$ exactly.
\begin{gather*}
C_{1,2M+1}^{\mathrm{odd}}(g) = \frac{1}{2}\sum_{n=-M}^{M} g_{2n}\\
C_{2,2M+1}^{\mathrm{odd}}(g) = \frac{1}{2}\sum_{l=-M}^{M} \sum_{|k|> M} g_{l+k} g_{l-k}
\end{gather*}
\end{lem}

\noindent\textbf{Proof.}
First of all, from \eqref{eq:C^odd} we have that
\begin{equation*}
C_{1,2M+1}^{\mathrm{odd}}(g) = \frac{1}{2} \int_{-\pi}^\pi g(x) S_{2M+1}(2x) \;\d x
\end{equation*}
\begin{multline*}
C_{2,2M+1}^{\mathrm{odd}}(g) = \frac{1}{2} \int_{-\pi}^\pi g^{2}(x) S_{2M+1}(2x) \;\d x\\
- \frac{1}{4}\int_{-\pi}^\pi \int_{-\pi}^\pi g(x) g(y) 2 S_{2M+1}(x+y) S_{2M+1}(x-y)\;\d x\;\d y
\end{multline*}
and using \eqref{eq:fourier expansion of S_(2M+1) kernel} we see that
\begin{equation*}
C_{1,2M+1}^{\mathrm{odd}}(g) = \frac{1}{2} \sum_{n=-M}^{M} g_{2n}
\end{equation*}
and
\begin{align*}
C_{2,2M+1}^{\mathrm{odd}}(g) &= \frac{1}{2}\sum_{l=-M}^{M} \sum_{k=-\infty}^\infty g_{k} g_{2l-k} - \frac{1}{2}\sum_{l=-M}^{M} \sum_{k=-M}^{M} g_{l+k} g_{l-k}\\
&= \frac{1}{2}\sum_{l=-M}^{M} \sum_{|k|> M} g_{l+k} g_{l-k}
\end{align*}
as required.
\qed

\begin{lem}\label{lem:SO(2M) odd}
For $\ell\geq 2$,
\begin{equation*}
\left|C_{\ell,2M+1}^{\mathrm{odd}}(g)\right| \leq \text{const}_\ell
\sum_{\substack{\mathbf{k}\in\Z^\ell\\|k_1|+\dots+|k_\ell|> 2M+1}}|g_{k_1}|\dots |g_{k_\ell}|
\end{equation*}
\end{lem}

\textbf{Proof.}
Fix $\sigma\in P(\ell,m)$, and for $\mathbf{k}=(k_1,\dots,k_\ell)\in\Z^\ell$ set
\begin{gather*}
K_1 = \sum_{l=1}^{\lambda_1} k_l\\
K_2 = \sum_{l=\lambda_1+1}^{\lambda_1+\lambda_2} k_l\\
\vdots\nonumber\\
K_m = \sum_{l=\lambda_1+\dots+\lambda_{m-1}+1}^{\ell} k_l
\end{gather*}
(recall that $\ell=\lambda_1+\dots+\lambda_m$).
Therefore
\begin{equation*}
\prod_{j=1}^m g^{\lambda_j}(x_j)
= \sum_{\mathbf{k}\in\Z^\ell} \prod_{l=1}^\ell g_{k_l} \prod_{j=1}^m e^{\i K_j x_j} 
\end{equation*}

Hence, the integral in \eqref{eq:C^odd}
\begin{multline*}
\int_{[-\pi,\pi]^m} \prod_{j=1}^m g^{\lambda_j}(x_j) S_{2M+1}(x_j-\epsilon_j x_{j+1}) \;\d x_j\\
= \sum_{-M \leq n_1,\dots,n_m \leq M} \sum_{\mathbf{k}\in \Z^\ell} \prod_{l=1}^\ell g_{k_l} \int_{[-\pi,\pi]^m} \prod_{j=1}^m  e^{\i K_j x_j} e^{\i n_j(x_j-\epsilon_j x_{j+1})} \frac{\d x_j}{2\pi}\\
= \sum_{\mathbf{k}\in \Z^\ell} \prod_{l=1}^\ell g_{k_l} \sum_{-M \leq n_1,\dots,n_m \leq M} \int_{[-\pi,\pi]^m} \prod_{j=1}^m   \exp\left(\i x_j (K_j+n_j-\epsilon_{j-1} n_{j-1})\right) \frac{\d x_j}{2\pi}
\end{multline*}
where we have used \eqref{eq:fourier expansion of S_(2M+1) kernel} to express $S_{2M+1}(x_j-\epsilon_j x_{j+1})$ in its Fourier representation, and we have defined $\epsilon_0=\epsilon_m$, $n_0=n_m$ (so all indices are cyclic).

The integral above will be $1$ or $0$ depending on whether $n_j-\epsilon_{j-1} n_{j-1}=-K_j$ or not, so defining
\begin{multline}\label{eq:defn_curly_N}
\mathcal{N}(M,\sigma,\mathbf{k},\mathbf{\epsilon}) = \\
\#\left\{-M \leq n_1,\dots,n_m \leq M \ : \ n_j-\epsilon_{j-1} n_{j-1} = -K_j \ , \  j=1,\dots,m \right\}
\end{multline}
(the $K_1,\dots,K_m$ depend on both $\mathbf{k}$ and $\sigma$, recall) we see that
\begin{equation}\label{eq:cumulant (odd) Fourier}
C_{\ell,2M+1}^{\mathrm{odd}}(g) =  \sum_{\mathbf{k}\in \Z^\ell} \prod_{l=1}^\ell g_{k_l} \sum_{m=1}^\ell \sum_{\sigma\in P(\ell,m)} \!\!\frac{(-1)^{m+1}(m-1)!}{2^m}\!\!\!\!\sum_{\substack{\epsilon_1,\dots,\epsilon_m=\pm 1\\\prod\epsilon_j=-1}}\!\!\! \mathcal{N}(M,\sigma,\mathbf{k},\mathbf{\epsilon})
\end{equation}

\begin{lem}\label{lem:Zeev's counting lemma}
Let $\prod_{j=1}^m \epsilon_j=-1$. Then $\mathcal{N}(M,\sigma,\mathbf{k},\mathbf{\epsilon})$ is either $0$ or $1$.
\begin{itemize}
\item If $\sum_{l=1}^\ell k_l$ is odd then $\mathcal{N}(M,\sigma,\mathbf{k},\mathbf{\epsilon})=0$.
\item If $\sum_{l=1}^\ell k_l$ is even and $\sum_{l=1}^\ell |k_l| \leq 2M$ then $\mathcal{N}(M,\sigma,\mathbf{k},\mathbf{\epsilon})=1$.
\end{itemize}
\end{lem}
(proof deferred until the end of this section).

Therefore, if $\sum_{l=1}^\ell |k_l| \leq 2M+1$ then
\begin{multline}\label{eq:these terms vanish}
\sum_{m=1}^\ell \sum_{\sigma\in P(\ell,m)} (-1)^{m+1}(m-1)!\frac{1}{2^m}\sum_{\substack{\epsilon_1,\dots,\epsilon_m=\pm 1\\\prod\epsilon_j=-1}} \mathcal{N}(M,\sigma,\mathbf{k},\mathbf{\epsilon}) \\
= \frac{1}{2}\mathcal{M}(\mathbf{k}) \sum_{m=1}^\ell \sum_{\sigma\in P(\ell,m)} (-1)^{m+1}(m-1)!
\end{multline}
where
\begin{equation*}
\mathcal{M}(\mathbf{k}) =
\begin{cases}
1 & \text{ if } \sum_{l=1}^\ell k_l \text{ is even}\\
0 & \text{ otherwise}
\end{cases}
\end{equation*}
Using the fact that for $\ell\geq 2$,
\begin{equation*}
\sum_{m=1}^\ell\sum_{\sigma\in P(\ell,m)} (-1)^{m+1}(m-1)!  = 0
\end{equation*}
we see that \eqref{eq:these terms vanish} vanishes for $\sum_{l=1}^\ell |k_l| \leq 2M+1$ if $\ell\geq 2$. Inserting this into \eqref{eq:cumulant (odd) Fourier} and estimating the contribution from the terms with $\sum_{l=1}^\ell |k_l| \geq 2M+2$ we see that
\begin{equation*}
\left|C_{\ell,2M+1}^{\mathrm{odd}}(g)\right| \leq \text{const}_\ell \sum_{\substack{\mathbf{k}\in\Z^\ell\\\sum_{l=1}^\ell |k_l| \geq 2M+2}} \prod_{l=1}^\ell |g_{k_l}|
\end{equation*}

This completes the proof of Lemma~\ref{lem:SO(2M) odd}.
\qed

\textbf{Proof of Lemma~\ref{lem:Zeev's counting lemma}}

We treat all indices as cyclic modulo $m$. So $n_0 = n_m$ and $n_{m+1}=n_1$ etc.

We assume that $\prod_{j=1}^m \epsilon_j = -1$.

Define the $m\times m$ matrix $E$ to be such that 
\begin{equation*}
E_{i,j} = 
\begin{cases}
\epsilon_{i-1} & \text{ if } j=i-1\\
0 & \text{ otherwise}
\end{cases}
\end{equation*}
so that
\begin{equation*}
(E \mathbf{n})_j = \epsilon_{j-1} n_{j-1}
\end{equation*}

From the definition of $\mathcal{N}(M,\sigma,\mathbf{k},\mathbf{\epsilon})$ (which is given in \eqref{eq:defn_curly_N}) we see that it is the number of solutions of $(I-E)\mathbf{n}=-\mathbf{K}$ subject to $-M\leq n_j\leq M$.

Now,
\begin{align*}
(E^k \mathbf{n})_j &= \epsilon_{j-1} (E^{k-1} \mathbf{n})_{j-1} \\
&= \epsilon_{j-1}\epsilon_{j-2}\dots\epsilon_{j-k} n_{j-k}
\end{align*}
and so $E^m = \epsilon_1\dots\epsilon_m I = -I$ by cyclicity of indices and the assumption that $\prod_{j=1}^m \epsilon_j = -1$.

Hence $2I = I-E^m$. But $I-E^m$ factorizes as
\begin{equation*}
I-E^m = (I-E)(I+E+\dots+E^{m-2}+E^{m-1})
\end{equation*}
and therefore
\begin{equation*}
(I-E)^{-1} = \frac{1}{2}(I+E+\dots+E^{m-2}+E^{m-1})
\end{equation*}

If we ignore the restriction that $-M\leq n_j\leq M$ then, over the reals, there is exactly one solution to $(I-E)\mathbf{n}=-\mathbf{K}$ which is
\begin{equation}\label{eq:soln for n_j}
n_{j} = -\frac{1}{2}\left(K_j+\epsilon_{j-1}K_{j-1}+\epsilon_{j-1}\epsilon_{j-2}K_{j-2}+\dots+\epsilon_{j-1}\epsilon_{j-2}\dots\epsilon_{j-m+1}K_{j-m+1}\right)
\end{equation}
This is a solution over the integers if $n_j$ is an integer, which will be the case when the term inside the bracket is even. Since $\epsilon_j\equiv 1(\mod 2)$ for all $j$, the term inside the bracket is even when
\begin{equation*}
K_j+K_{j-1}+\dots+K_{j-m+1} = \sum_{i=1}^m K_i = \sum_{l=1}^\ell k_l
\end{equation*}
is even. There are no solutions over the integers when this is odd. (Note that the even and oddness is independent of $\mathbf{\epsilon}$ and of the partition $\sigma$).

Finally, one must check that the condition $-M\leq n_j \leq M$ holds. From \eqref{eq:soln for n_j} we see that
\begin{equation*}
|n_j| \leq \frac{1}{2} \sum_{i=1}^m |K_i| \leq \frac{1}{2}\sum_{l=1}^\ell |k_l|
\end{equation*}
and so if we assume that $\sum_{l=1}^\ell |k_l| \leq 2M$, then the condition holds.

Thus $\mathcal{N}(M,\sigma,\mathbf{k},\mathbf{\epsilon})=0$ if $\sum_{l=1}^\ell k_l$ is odd, and $\mathcal{N}(M,\sigma,\mathbf{k},\mathbf{\epsilon})=1$ if $\sum_{l=1}^\ell k_l$ is even and $\sum_{l=1}^\ell |k_l| \leq 2M$.

This proves Lemma~\ref{lem:Zeev's counting lemma}.
\qed

\subsection{The calculation of $C_{\ell,2M}^{\mathrm{odd}}(g)$}

Basically, this section is like the previous, with the essential change being that
$$S_{2M}(z) = \frac{1}{2\pi} \sum_{\substack{n=-(2M-1)\\n \text{ odd}}}^{2M-1} e^{-\i n z/2}$$
as opposed to \eqref{eq:fourier expansion of S_(2M+1) kernel} which says
$$S_{2M+1}(z) = \frac{1}{2\pi} \sum_{\substack{n=-2M\\n \text{ even}}}^{2M} e^{-\i n z/2}$$

\begin{lem}
One can calculate $C_{1,2M}^{\mathrm{odd}}(g)$ and $C_{2,2M}^{\mathrm{odd}}(g)$ exactly.
\begin{equation*}
C_{1,2M}^{\mathrm{odd}}(g) = \frac{1}{2}\sum_{n=-(M-1)}^{M} g_{2n-1}
\end{equation*}
\begin{equation*}
C_{2,2M}^{\mathrm{odd}}(g) = \frac{1}{2}\sum_{\substack{n=-(2M-1)\\n \text{ odd}}}^{2M-1} \sum_{\substack{|m|\geq 2M+1\\m \text{ odd}}} g_{\frac{1}{2}(n+m)} g_{\frac{1}{2}(n-m)}
\end{equation*}
\end{lem}

\begin{lem}\label{lem:SO(2M+1) odd}
For $\ell\geq 2$,
\begin{equation*}
\left|C_{\ell,2M}^{\mathrm{odd}}(g)\right| \leq \text{const}_\ell
\sum_{\substack{\mathbf{k}\in\Z^\ell\\|k_1|+\dots+|k_\ell|> 2M}}|g_{k_1}|\dots |g_{k_\ell}|
\end{equation*}
\end{lem}

The proof goes through the same as before, with equation \eqref{eq:defn_curly_N} becoming
\begin{multline*}
\mathcal{N}_{odd}(M,\sigma,\mathbf{k},\mathbf{\epsilon}) = \#\bigl\{-(2M-1) \leq n_j \leq 2M-1 \ , \ n_j \text{ odd}\ : \\
 \tfrac{1}{2}n_j-\epsilon_{j-1}\tfrac{1}{2} n_{j-1} = -K_j \ , \  j=1,\dots,m \bigr\}
\end{multline*}
Rewriting equation \eqref{eq:soln for n_j} we see the solution requested by $\mathcal{N}_{odd}(M,\sigma,\mathbf{k},\mathbf{\epsilon})$ is
\begin{equation*}
n_{j} = -\left(K_j+\epsilon_{j-1}K_{j-1}+\epsilon_{j-1}\epsilon_{j-2}K_{j-2}+\dots+\epsilon_{j-1}\epsilon_{j-2}\dots\epsilon_{j-m+1}K_{j-m+1}\right)
\end{equation*}
so long as $n_j$ is odd and $-(2M-1)\leq n_j \leq 2M-1$ (and there is no solution otherwise). Therefore Lemma~\ref{lem:Zeev's counting lemma} becomes
\begin{lem}\label{lem:counting lemma for odd orthog}
Let $\prod_{j=1}^m \epsilon_j=-1$. Then $\mathcal{N}_{odd}(M,\sigma,\mathbf{k},\mathbf{\epsilon})$ is either $0$ or $1$.
\begin{itemize}
\item If $\sum_{l=1}^\ell k_l$ is even then $\mathcal{N}_{odd}(M,\sigma,\mathbf{k},\mathbf{\epsilon})=0$.
\item If $\sum_{l=1}^\ell k_l$ is odd and $\sum_{l=1}^\ell |k_l| \leq 2M-1$ then $\mathcal{N}_{odd}(M,\sigma,\mathbf{k},\mathbf{\epsilon})=1$.
\end{itemize}
\end{lem}

\section{Moments of traces}

We will now use Theorem~\ref{thm:cumulants SO(N)} to prove the second part of Theorem~\ref{thm:mmts traces}. (The proof of the first part from Theorem~\ref{thm:cumulants Sp(2M)} being analogous).

Recall from \eqref{eq:defn mu} that
\begin{equation*}
C_{\ell}^{\SO(N)}(g) = \sum_{n\in\Z^\ell} \mu_\ell^{\SO(N)}(n_1,\dots,n_\ell) \prod_{j=1}^\ell g_{n_j}
\end{equation*}
where $\mu_\ell^{\SO(N)}(n_1,\dots,n_\ell)$ is invariant under permutations of its arguments. Assuming $g_0=0$, then we have
\begin{itemize}
\item If $|n_1|<N$ then
\begin{equation*}
\mu_1^{\SO(N)}(n_1) = 
\begin{cases}
1 & \text{ if $n_1\neq 0$ is even}\\
0 & \text{ otherwise}
\end{cases}
\end{equation*}
\item If $|n_1|+|n_2| < N$ then
\begin{equation*}
\mu_2^{\SO(N)}(n_1,n_2) = 
\begin{cases}
|n_1| & \text{ if $|n_1|=|n_2|$}\\
0 & \text{ otherwise}
\end{cases}
\end{equation*}
\item If  $\ell\geq 3$ and $\sum_{j=1}^\ell |n_j| < N$ then $\mu_{\ell}^{\SO(N)}(n_1,\dots,n_\ell)=0$.
\end{itemize}

It is also true that if $g_0=0$,
\begin{align}
\E_{G}\bigl\{ &(\Tr g(U) -C_1^{G}(g))^m\bigr\} \nonumber\\
&= 2^m \sum_{n\in\N^m} \E_{G}\left\{(\Tr U^{n_1}-\mu_1^G(n_1)) \dots (\Tr U^{n_m}-\mu_1^G(n_m))\right\} \prod_{j=1}^m g_{n_j}\label{eq:moments written out}\\
&= \sum \left(\frac{C_2^G(g)}{2!}\right)^{k_2} \left(\frac{C_3^G(g)}{3!}\right)^{k_3}\dots \left(\frac{C_m^G(g)}{m!}\right)^{k_m} \frac{m!}{k_2! k_3! \dots k_m!} \label{eq:moments in terms of cumulants}
\end{align}
where the second sum runs over all values of $k_j\geq 0$ such that $\sum_{j=2}^m j k_j = m$ (it is simply writing the $m$th moment in terms of its cumulants, having subtracted the mean).

Let $a_j \in \{0,1,2\dots\}$ for $j=1,2,\dots$ by such that $\sum j a_j < N$. Define
\begin{equation*}
\eta_j=
\begin{cases}
1 & \text{ for even $j$}\\
0 & \text{ for odd $j$}
\end{cases}
\end{equation*}
so that $\mu_1^{\SO(N)}(j)=\eta_j$ for $|j|<N$.

Putting $m=\sum a_j$, we will evaluate the  coefficient of $\prod (g_{j})^{a_j}$ in \eqref{eq:moments written out} and in \eqref{eq:moments in terms of cumulants}, the two being equal to each other.

Consider first equation \eqref{eq:moments written out}. The coefficient of $\prod (g_{j})^{a_j}$ in
\begin{equation*}
2^m \sum_{n\in\N^m} \E_{\SO(N)}\left\{(\Tr U^{n_1}-\eta_{n_1}) \dots (\Tr U^{n_m}-\eta_{n_m})\right\} \prod_{j=1}^m g_{n_j}
\end{equation*}
equals 
\begin{equation}\label{eq:coeff on LHS}
\frac{2^m m!}{\prod (a_j)!} \E_{\SO(N)}\left\{ \prod (\Tr U^j-\eta_{j})^{a_j}\right\}
\end{equation}

Consider next equation \eqref{eq:moments in terms of cumulants}. Note that the restriction on the $a_j$ means that there is no contribution to the coefficient of $\prod (g_{j})^{a_j}$ from $C_{\ell}^{\SO(N)}(g)$ for all $\ell\geq 3$. Therefore the coefficient in \eqref{eq:moments in terms of cumulants} is $0$ if $m$ is odd and is the coefficient of $\prod (g_{j})^{a_j}$ in
\begin{equation*}
\frac{m!}{2^{m/2} (m/2)!} \left(C_2^{\SO(N)}(g)\right)^{m/2} = \frac{m!}{2^{m/2} (m/2)!} 2^m \sum_{n\in\N^m} \prod_{j=1}^{m/2} \mu_2^{\SO(N)}(n_{2j-1},n_{2j}) \prod_{j=1}^m g_{n_j}
\end{equation*}
if $m$ is even.
This coefficient is zero unless all the $a_j$ are even, in which case it is
\begin{equation}\label{eq:coeff on RHS}
\frac{m!}{(m/2)!} 2^{m/2} \frac{(m/2)!}{\prod (a_j/2)!} \prod j^{a_j/2}
\end{equation}
(to see this, note that the structure of $\mu_2^{\SO(N)}$ means that $n_{2j}$ must equal $n_{2j-1}$ for $j=1,\dots,m/2$. The second prefactor is just the number of ways of picking $m/2$ integers such that $a_1/2$ of them equal $1$, $a_2/2$ of them equal $2$ etc.).

Setting \eqref{eq:coeff on LHS}=\eqref{eq:coeff on RHS} and recalling that $m=\sum a_j$, we have
\begin{align*}
 \E_{\SO(N)}\left\{ \prod (\Tr U^j-\eta_{j})^{a_j}\right\} &= 
\begin{cases}
\prod j^{a_j/2} \frac{(a_j)!}{2^{a_j/2} (a_j/2)!} & \text{ if all the $a_j$ are even}\\
0 & \text{ otherwise}
\end{cases}\\
&= E\left\{ \prod (\sqrt{j} Z_j)^{a_j} \right\}
\end{align*}
where $Z_j$ are iid normal random variables with mean $0$ and variance $1$.

Observe that this can all be rewritten as
\begin{equation*}
 \E_{\SO(N)}\left\{ \prod (\Tr U^j)^{a_j}\right\} = E\left\{ \prod (\sqrt{j} Z_j + \eta_j)^{a_j} \right\}
\end{equation*}
and is valid so long as $\sum ja_j < N$.

\end{document}